\documentclass{amsart}
		\title[On the Bateman-Horm Conjecture
 about Polynomial Rings]{On the Bateman-Horm Conjecture
\\ about Polynomial Rings}
		\author[Bary-Soroker and Jarden]{Lior Bary-Soroker$^*$
and Moshe Jarden$^{**}$
}
		\address{Lior Bary-Soroker\\
		School of Mathematics\\
               	Tel Aviv University\\
Ramat Aviv, Tel Aviv 69978, Israel}
                \email{barylior@post.tau.ac.il}
      		\urladdr{http://www.tau.ac.il/~barylior/}

		\address{Moshe Jarden\\
		School of Mathematics\\
               	Tel Aviv University\\
Ramat Aviv, Tel Aviv 69978, Israel}
                \email{jarden@post.tau.ac.il}
      		\urladdr{http://www.tau.ac.il/~jarden/}
		
\thanks{$^*$ Research supported by the Alexander von Humboldt Foundation.}
\thanks{$^{**}$ Research supported by the Minkowski Center for
Geometry at Tel Aviv University, established by the Minerva
Foundation, and by an ISF grant.}

\input xy
\xyoption{all}
\usepackage{amssymb}

\def\phi{\varphi}
\def\alp{\alpha}
\def\bet{\beta}
\def\bbA{\mathbb{A}}
\def\bbF{\mathbb{F}}
\def\bbP{\mathbb{P}}
\def\bbQgal{{\tilde Q}}
\def\bbR{\mathbb{R}}
\def\bbZ{\mathbb{Z}}
\def\bfA{{\bf A}}
\def\bfa{{\bf a}}
\def\bfb{{\bf b}}
\def\bfbhat{{\hat\bfb}}
\def\bfo{{\bf o}}
\def\bfp{{\bf p}}
\def\bfx{{\bf x}}
\def\bhat{{\hat b}}
\def\cf{}
\def\chr{{\rm char}}
\def\colon{{:}\;}
\def\congr{\equiv}
\def\cont{\subseteq}
\def\Del{\Delta}
\def\del{\delta}
\def\eps{\epsilon}
\def\Fhat{\hat F}
\def\frp{\mathfrak{p}}
\def\Gal{{\rm Gal}}
\def\Gam{\Gamma}
\def\ghat{{\hat g}}

\def\hhat{{\hat h}}
\def\irr{{\rm irr}}
\def\isom{\cong}
\def\Ker{{\rm Ker}}
\def\Kgag{{\bar K}}
\def\Kgal{{\tilde K}}
\def\Lam{\Lambda}
\def\lam{\lambda}
\def\notdivide{\nmid}
\def\ome{\omega}
\def\onto{\mapsto}
\def\phihat{\hat\phi}
\def\Quot{{\rm Quot}}
\def\res{{\rm res}}
\def\Shat{{\hat S}}
\def\sig{\sigma}
\def\Spec{{\rm Spec}}
\def\st{\mathop{|\;}}
\def\Tet{\Theta}
\def\tet{\theta}
\def\tetgag{\bar\tet}
\def\tethat{\hat\tet}
\def\that{{\hat t}}
\def\total.deg{{\rm total.deg}}
\def\union{\bigcup}
\def\What{{\hat W}}
\def\xgag{{\bar x}}

\newbox\jtembox
\newdimen\jtemindent \jtemindent=25pt
\newdimen\oldindent  \oldindent=25pt
\def\jtem#1{\oldindent=\parindent\parindent=\jtemindent%
            \par\hangindent\parindent\indent\llap{\rm #1\enspace}%
            \parindent=\oldindent\ignorespaces}
\def\condition(#1){%
   \setbox\jtembox=\hbox{(#1)\enspace}\jtemindent=\wd\jtembox%
   \jtem{\rm(#1)}%
}

\def\eqalign#1{\null\,\vcenter{\openup1\jot
\mathsurround=0pt
\ialign{\strut\hfil$\displaystyle{##}$&$\displaystyle{}##$\hfil
\crcr#1\crcr}}\,}
\def\demo#1{\medskip\noindent{\rm #1.}\enspace}
\def\hefresh{
   \def\hefreshD{\mathop{\raise1.5pt\hbox{${\smallsetminus}$}}}
   \def\hefreshS{\mathop{\raise0.85pt\hbox{$\scriptstyle\smallsetminus$}}}
   \mathchoice{\hefreshD}{\hefreshD}{\hefreshS}{\hefreshS}}
\newif\ifextra \extratrue

\def\proclaimm#1#2#3{\ifdim\lastskip<\medskipamount\removelastskip
\medskip\fi\penalty-55\noindent#1#2#3\ignorespaces}
\def\proclaim #1{\proclaimm{\bf #1.\enspace}{}\sl}

\def\endproclaim{\par\rm
 \ifdim\lastskip<\smallskipamount\removelastskip\penalty55\smallskip\fi}
\def\rem#1{\proclaimm{\sl#1:}{}{\quad\rm}}
\def\reminfo#1#2{\proclaimm{\sl#1.}{\enspace\sl#2.}{\quad\rm}}
\def\endrem{\qed}
\newbox\keybox
\newdimen\keywidth
\def\references#1{
\setbox\keybox=\hbox{[#1]\enspace}\keywidth=\wd\keybox
\oldindent=\parindent \parindent=\keywidth
\frenchspacing\rm
   \ifdim\lastskip<\bigskipamount\removelastskip\bigskip\fi
\centerline{\bf References}\medskip\nobreak}
\def\ref[#1]{\par\smallskip\hangindent\parindent\indent\llap{\hbox to\keywidth%
    {[#1]\hfil\enspace}}\ignorespaces}

\def\subdemo#1{\proclaimm{\sc#1\unskip\ifextra:\fi\extratrue}{}%
{\quad\rm}}
\def\subdemoinfo#1#2{\proclaimm{\sc#1\ifextra:\fi\extratrue}%
{\enspace\sl #2\ifextra.\fi\extratrue}{\quad\rm}}

\theoremstyle{plain}

\begin{document}
 
\begin{abstract}
Given a power $q$ of a prime number $p$ and ``nice''
polynomials $f_1,\ldots,f_r\in\bbF_q[T,X]$
with $r=1$ if $p=2$, we establish an
asymptotic formula for the number of pairs
$(a_1,a_2)\in\bbF_q^2$ such that
$f_1(T,a_1T+a_2),\ldots,f_r(T,a_1T+a_2)$
are irreducible in $\bbF_q[T]$.
In particular that number tends to infinity
with $q$.
\end{abstract}

\maketitle
\noindent
{\bf Introduction}

Let $f_1,\ldots,f_r\in\bbZ[X]$ be non-associate irreducible
polynomials with positive
leading coefficients.
A conjecture of Bateman and Horn [BaH62, (1)]
predicts for $x>1$ that the number
$N(f_1,\ldots,f_r;x)$ of positive integers $1\le n\le x$
such that $f_1(n),\ldots,f_r(n)$ are prime
numbers satisfies
$$
N(f_1,\ldots,f_r;x)\sim
{s(f_1,\ldots,f_r)\over\prod_{i=1}^r\deg(f_i)}{x\over\log^rx},
$$
where
$$
\eqalign{
s(f_1,\ldots,f_r)
&=\prod_p{1-{\ome(p)\over p}
\over\big(1-{1\over p}\big)^r}
\hbox{ and }
\cr
\ome(p)
&=\#\{0\le n\le p-1\st f_1(n)\cdots f_r(n)\congr0\mod p\}.
}
$$
If $\ome(p)=p$ for some $p$, then $s(f_1,\ldots,f_r)=0$.
Also, for each $n\in\bbZ$ there exists $1\le i\le r$ such
that $p|f_i(n)$, thus $N(f_1,\ldots,f_r;x)$ is bounded.
If $\ome(p)<p$ for all $p$, then $s(f_1,\ldots,f_r)$
converges to a positive real number [BaH62, p.~364],
hence the conjecture
predicts the existence of
infinitely many $n$'s such that $f_1(n),\ldots,f_r(n)$ are all
prime numbers.
This is a conjecture of Schinzel.
The special case where $r=1$, $f_1(X)=X$, and $s(f_1)=1$
reduces to the prime number theorem.
When $r=1$ and $f_1(X)=aX+b$ with $\gcd(a,b)=1$, and
$s(f_1)={a\over \phi(a)}$, we get Dirichlet's theorem.
Likewise, in the case where $f_1(X)=X$
and $f_2(X)=X+2$, the Bateman-Horn Conjecture generalizes the
Hardy-Littlewood conjecture about the density of the twin
primes.
Calculations made by Littlewood show that in this case
$s(f_1,f_2)\approx1.32$ [BaT04, p.~335]. 

It is customary in number theory to replace $\bbZ$ by the
ring $\bbF_q[T]$, where $q$ is a power of a prime number
$p$.
In our case we would like to consider irreducible
non-associate polynomials $f_1,\ldots,f_r\in\bbF_q[T,X]$ and
ask for the number of $g\in\bbF_q[T]$ such that
$f_1(T,g(T)),\ldots,f_r(T,g(T))$ are irreducible in
$\bbF_q[T]$.

It turns out that the naive analog of Schinzel's conjecture
may fail when one of the $f_i$'s is a polynomial in $X^p$.
For example, Swan proves in [Swa62] that $g(T)^8+T^3$ is
reducible for each $g\in\bbF_2[T]$.
Thus, more restrictions are needed to restore Schinzel's
and the Bateman-Horn conjectures [CCG08].

In a conference in the American Institute of Mathematics
[Gao03], S. Gao posed the following question:

\proclaim {Problem A}
Let $f\in\bbF_q[T,X]$ be an irreducible polynomial.
Count (or estimate) the number of pairs
$(a_1,a_2)\in\bbF_q\times\bbF_q$ such that the polynomial
$f(T,a_1T+a_2)$ is irreducible in $\bbF_q[T]$.
\endproclaim

Bender and Wittenberg [BeW05, Thm.~1.1 and Prop.~4.1]
prove the first result in this direction:

\proclaim {Proposition B}
Let $q$ be a power of a prime number and $f_1,\ldots,f_r\in \bbF_q[T,X]$
polynomials of degrees $d_1,\ldots,d_r$, respectively.
Suppose each $1\le i\le r$ satisfies
\condition (\cf1a)
$p\notdivide d_i(d_i-1)$ and
\condition (\cf1b)
the Zariski closure $C_i$ in $\bbP_{\bbF_q}^2$ of the affine plane curve
defined by $f_i(T,X)=0$ is smooth.

Then for each large $k$ there exists
$(a_1,a_2)\in\bbF_{q^k}^2$ such that each of the polynomials
$f_1(T,a_1T+a_2),\ldots,f_r(T,a_1T+a_2)$ is irreducible in
$\bbF_{q^k}[T]$.

Moreover, for $r=1$ and $q>9(d(d-1)d+2)^2$ the number of
pairs $(a_1,a_2)\in\bbF_q^2$ such that $f_1(T,a_1T+a_2)$ is
irreducible in $\bbF_q[T]$ is at least
${1\over d!}q^2+c_1q^{3/2}+c_2q+c_3q^{1\over2}+c_4$,
where $c_1,c_2,c_3,c_4$ are explicitely given constants
depending only on $d$.
\endproclaim

Note that Condition (\cf1a) implies that $p\ne2$.
Our main
result improves Proposition~B in three ways.
First, it includes the
case $p=2$, albeit with the restriction that $r=1$ in this case.
Second, we
replace the condition that the $C_i$'s be smooth by a less
restrictive condition of being ``characteristic-$0$-like''
(see Section 1 for definition)
and {\bf nodal} (i.e.~having only nodes as singularities).
Finally, our result is quantitative for arbitrary $r$, that is, we
give an asymptotic formula for the number of pairs $(a_1,a_2)$ for
which the polynomials $f_i(T,a_1T+a_2)$ are irreducible,
when $q\to\infty$.
Thus, our result solves Problem A for ``nice'' polynomials.

\proclaim {Theorem C}
Let
$f_1,\ldots,f_r\in\bbF_q[T,X]$ be 
non-associate absolutely irreducible characteristic-$0$-like
nodal polynomials, where
$q$ is a power of a prime $p$ such that $r=1$ if $p=2$.
For each $1\le i\le r$ let $d_i$ be the degree of $f_i$
and set $d=\prod_{i=1}^rd_i$.
Then,
$$
\eqalign{
\#\{(a_1,a_2)\in\bbF_q^2&\st f_i(T,a_1T+a_2)
\hbox{ is irreducible in } \bbF_q[T],
\ i=1,\ldots,r\}
\cr
&={q^2\over d}+O(q^{3/2}),
}
$$
where the constant of the $O$ is a computable function in
$\sum_{i=1}^r d_i$.
\endproclaim

Note that the leading terrm in our approximation formula is
$q^2\over d$.
For large $q$ and for $d>2$, this improves the lower bound
given in Theorem B for $r=1$ that has 
the leading term $q^2\over d!$.

Apart from standard combinatorial arguments, the proof of
Theorem C is based on three ingredients:
geometrical arguments used in the proof of the theorem about
the stability of fields [FrJ08, Section 18.9], the field
crossing argument [FrJ08, Section 24.1], and the Lang-Weil
estimates.

Proposition 4.1 proves that Proposition B is a
special case of Theorem C.

Finally we mention that Pollack [Pol08] and the first author [BaS10]
treat the case of the analog of
the Bateman--Horn Conjecture  when $T$ does
not occur in the $f_i$'s.


\section{Direct Product of Symmetric Groups}

Let $\Gam$ be a projective plane curve defined over a field
$K$ by an absolutely irreducible homogeneous equation
$f(X_0,X_1,X_2)=0$
with a generic point $\bfx=(x_0{:}x_1{:}x_2)$.
Then the point
$$
\bfx^*=(x_0^*{:}x_1^*{:}x_2^*)
=\Big({\partial f\over\partial X_0}(\bfx)
:{\partial f\over\partial X_1}(\bfx)
:{\partial f\over\partial X_2}(\bfx)\Big)
\leqno (\cf1)
$$
is a generic point of an absolutely irreducible projective
plane $K$-curve $\Gam^*$, known as the {\bf dual curve} of
$\Gam$.
The points of $\Gam^*$ parametrize the tangents of $\Gam$ at
simple points.
The map $\bfx\onto\bfx^*$ extends to a rational map
$\Gam\to\Gam^*$.

\proclaim {Lemma 1.1}
Let $K$ be an algebraically closed field with $\chr(K)\ne2$.
Let $\Gam$ and $\Del$ be distinct 
absolutely irreducible projective plane $K$-curves.
Suppose both $\Gam$ and $\Del$ have only finitely many
inflection points.
Then $\Gam$ and $\Del$ have only finitely many common
tangents.
\endproclaim

\demo {Proof}
Assume $\Gam$ and $\Del$ have infinitely many common
tangents.
Then $\Gam^*(K)\cap\Del^*(K)$ is infinite, hence
$\Gam^*=\Del^*$, so $\Gam^{**}=\Del^{**}$.
Since both $\Gam$ and $\Del$ have only finitely many
inflection points, $\Gam^{**}=\Gam$ and $\Del^{**}=\Del$
[GeJ89, Prop.~4.5].
Therefore, $\Gam=\Del$, in contrast to our assumption.
\qed

Given an absolutely irreducible projective plane
$K$-curve $\Gam$, we write
$\Gam_\Kgal$ for the
$\Kgal$-curve obtained by base change from $K$ to
its algebraic closure $\Kgal$.
It is well known, that if $\Gam$ is not a line and $\chr(K)=0$, then
\condition (\cf2a)
$\Gam_\Kgal$ has only finitely many inflection points,
\condition (\cf2b)
$\Gam_\Kgal$ has only finitely many double tangents, and
\condition (\cf2c)
$\Gam_\Kgal$ is not strange,

\noindent
[FrJ76, Lemma 3.2].
In general we
say that $\Gam$ is a {\bf characteristic-$0$-like curve}
if it satisfies Condition (\cf2).

We say that $\Gam$ is a {\bf nodal curve}, if all of the
singular points of $\Gam_\Kgal$ are nodes [Ful89, p.~66].

We say that an absolutely irreducible polynomial $f\in
K[T,X]$ is {\bf charac\-teristic-$0$-like nodal}
if the Zariski closure $\Gam$ in
$\bbP_K^2$ of the affine plane curve defined by the equation
$f(T,X)=0$ is a characteristic-$0$-like nodal curve.

Finally we say that polynomials $f,g\in K[T,X]$ are {\bf
non-associate} if $f\ne\lam g$ for all $\lam\in K^\times$.

\proclaim {Lemma 1.2}
Let $K$ be an infinite field and $f_1,\ldots,f_r\in K[T,X]$
 absolutely irreducible non-associate
characteristic-$0$-like nodal polynomials of degrees
$d_1,\ldots,d_r$, respectively.
Then:
\condition (a)
There exist $\alp,\bet\in K$ such that for each $1\le i\le
r$ we have
$$
\Gal\Big(f_i\big(T,{\alp+T+\bet U\over
U}\big),\Kgal(U)\Big)\isom S_{d_i}.
$$
\condition (b)
If, in addition, $\chr(K)\ne2$ and we set $f=f_1\cdots f_r$,
then $\alp,\bet$ can be chosen such that
$$
\Gal\Big(f\big(T,{\alp+T+\bet U\over
U}\big),\Kgal(U)\Big)\isom\prod_{i=1}^rS_{d_i}.
$$
\endproclaim

\demo {Proof of (a)}
For each $1\le i\le r$
let $\Gam_i$ be the Zariski closure in $\bbP_K^2$ of the
absolutely irreducible affine plane $K$-curve defined over
$K$ by the equation $f_i(T,X)=0$.
By our assumption on the $f_i$'s, the absolutely irreducible
projective plane $K$-curves
$\Gam_1,\ldots,\Gam_r$ are distinct.
Since the $\Gam_i$'s are characteristic-$0$-like curves,
there are finitely many lines
$L_1,\ldots,L_m$ in $\bbP_\Kgal$ such that if
$$
\bfo\in\bbP^2(\Kgal)\hefresh\union_{i=1}^r\Gam_i(\Kgal)\cup
\union_{j=1}^mL_j(\Kgal),
\leqno (\cf3)
$$
then for each $1\le i\le r$ we have [FrJ76, proof of Lemma
3.2]:
\condition (\cf4a)
$\bfo$ lies on no tangent to $\Gam_i$ at an inflection
point;
\condition (\cf4b)
no double tangent to $\Gam_i$ goes through $\bfo$;
\condition (\cf4c)
$\bfo$ lies on no line that goes through two singular points
of $\Gam_i$;
\condition (\cf4d)
$\bfo$ lies on no tangent to $\Gam_i$ that goes through a
singular point of $\Gam_i$; and
\condition (\cf4e)
only finitely many lines through $\bfo$ are tangents to
$\Gam_i$.

Since $K$ is infinite, we may choose a point
$\bfo=(1{:}{-}\alp{:}\bet)\in\bbP^2(K)$ that satisfies (\cf3),
hence also (\cf4).
Then $\bfo$ is the intersection of the lines $\Lam_\alp$ and
$\Lam_\bet$ respectively defined by the homogeneous
equations $-\alp X_0-X_1=0$ and $\bet X_0-X_2=0$ with
coefficients in $K$.

We consider $1\le i\le r$.
Then $\deg(\Gam_i)=\deg(f_i)=d_i$.
We choose a transcendental element $t$ over $K$ and a root
$x_i$ of the equation $f_i(t,X)$ in $\widetilde{K(t)}$,
and set $F_i=K(t,x_i)$.
The projection
$\lam\colon\bbP_\Kgal^2\hefresh\{\bfo\}\to\bbP_K^1$ from
$\bfo$ is defined over $K$ by
$$
\lam(X_0{:}X_1{:}X_2)=(-\alp X_0-X_1:\bet X_0-X_2).
$$
In particular, $\lam(1{:}t{:}x_i)=(-\alp-t:\bet-x_i)$ is a
generic point of $\bbP_K^1$ and
$$
u_i={-\alp-t\over\bet-x_i}
\leqno (\cf5)
$$
is transcendental over $K$ (because $\deg(f_i)\ge2$).

A central argument in the proof of [FrJ76, Lemma 3.3]
states that every line in $\bbP_\Kgal^2$
that passes through $\bfo$ cuts $\Gam_{i,\Kgal}$
in at least $d_i-1$ points and almost every such line
cuts $\Gam_{i,\Kgal}$ in $d_i$ points.
Then, by [FrJ76, Lemma 2.1], $u_i$ is a
separating transcendental element for $F_i/K$ and the Galois
closure $\Fhat_i$ of $F_i/K(u_i)$ satisfies
$$
\Gal(\Fhat_i/K(u_i))\isom\Gal(\Fhat_i\Kgal/\Kgal(u_i))\isom
S_{d_i}.
$$
By (\cf5), 
$x_i={\alp+t+\bet u_i\over u_i}$, so $F_i=K(t,u_i)$.
Since $f_i(t,x_i)=0$, we have
$$
f_i\big(t,{\alp+t+\bet u_i\over u_i}\big)=0.
$$
Hence, the rational function
$
f_i\big(T,{\alp+T+\bet U\over U}\big)
$
in the variables $T,U$
is absolutely irreducible, separable in $T$, and with Galois
group over $\Kgal(U)$ isomorphic to $S_{d_i}$, as claimed.

\demo {Proof of (b)}
Now we assume that $\chr(K)\ne2$.
Then, by Lemma 1.1 only finitely many lines in
$\bbP_\Kgal^2$ are tangents to two of the curves
$\Gam_1,\ldots,\Gam_r$.
Using the notation of the proof of (a), we may assume
that these lines belong to the
set $\{L_1,\ldots,L_m\}$.
Then the point $\bfo$ satisfies, in addition to (\cf4), also
the following condition:
\condition (\cf6)
Each line through $\bfo$ is a tangent to at most one of the
curves $\Gam_1,\ldots,\Gam_r$.

Let $1\le i\le r$ and consider a prime divisor $\frp$ of
$\Kgal(u_i)/\Kgal$ that ramifies in $\Fhat_i\Kgal$.
Then $\frp$ may be identified with the intersection point of
$\bbP_\Kgal^1$ with
a tangent $L$ to $\Gam_i$ at a point $P$ that goes through $\bfo$.
Indeed, in this case, the intersection multiplicity of $L$
with $\Gam_i$
at $P$ is $2$ and is $1$ at all other intersection points
(by (\cf4)).
Thus, by Bezout's theorem, $L$ has $d_i-1$ intersection
points with $\Gam_i$, so $\frp$ decomposes in $F_i\Kgal$ as
$\frp=2\frp_1+\frp_2+\cdots+\frp_{d-1}$ with distinct prime divisors
$\frp_1,\cdots,\frp_{d-1}$.
We also identify $u_i$ with the variable $U$.
It follows from (\cf6) that for all $i\ne j$,
the set of the prime divisors of $\Kgal(U)/\Kgal$ that
ramify in $\Fhat_i\Kgal$ is disjoint from the set of prime
divisors of $\Kgal(U)/\Kgal$ that ramify in $\Fhat_j\Kgal$.
Therefore, by Riemann-Hurwitz, we have for each $1\le j\le r$ that
$(\Fhat_j\Kgal)\cap(\prod_{i\ne j}\Fhat_j\Kgal)=\Kgal(U)$
[FrJ08, Rem.~16.2(b)].
Consequently, by (a),
$$
\Gal\Big(f\big(T,{\alp+T+\bet U\over U}\big),\Kgal(U)\Big)
\isom\prod_{i=1}^r\Gal(\Fhat_i/\Kgal(U))
\isom\prod_{i=1}^rS_{d_i},
$$
as claimed.
\qed

\smallskip
The following corollary does not assume $K$ to be infinite.

\proclaim {Corollary 1.3}
Let $K$ be a field and let
$f_1,\ldots,f_r\in K[T,X]$ be non-associate
characteristic-$0$-like nodal
absolutely irreducible
polynomials of degrees $d_1,\ldots,d_r$,
respectively.
Let $A$ and $B$ be variables.
We set $f=f_1\cdots f_r$ and assume that $r=1$ if
$\chr(K)=2$.
Then,
$\Gal(f(T,AT+B),K(A,B))\isom\prod_{i=1}^rS_{d_i}$.
\endproclaim

\demo {Proof}
First we prove the corollary for $\Kgal$ rather than for $K$.
To this end we use the fact that $\Kgal$ is
infinite to choose $\alp,\bet\in\Kgal$ such that 
(a) of Lemma
1.2 holds,
and also (b) of that Lemma holds
if $\chr(K)\ne2$.
Then we extend the specialization
$(A,B)\to\big({1\over U},{\alp+\bet U\over U}\big)$
to a
$\Kgal$-place $\phi$ of the splitting field of $f(T,AT+B)$
over $K(A,B)$.
Then
$\Gal\big(f\big(T,{\alp+T+\bet U\over U}\big),\Kgal(U)\big)$ is a
quotient of a subgroup of
$\Gal(f(T,AT+B),\Kgal(A,B))$, namely the quotient
of the decomposition group by the inertia group of $\phi$.
Since, $\Gal(f(T,AT+B),\Kgal(A,B))\le\prod_{i=1}^rS_{d_i}$
and 
$\Gal\big(f(T,{\alp+T+\bet U\over U}),\Kgal(U)\big)
\isom\prod_{i=1}^r S_{d_i}$,
we conclude that $\Gal(f(T,AT+B),\Kgal(A,B))\isom\prod_{i=1}^rS_{d_i}$.

The truth of the corollary for $K$ now follows from its truth for
$\Kgal$ and from the
following inclusion of groups:
$
\prod S_{d_i}
\isom\Gal(f(T,AT+B),\Kgal(A,B))
\le\Gal(f(T,AT+B),K(A,B))
\le\prod_{i=1}^rS_{d_i}.
$
\qed

\section{The Field Crossing Argument}

The field crossing argument has already been used in the
original proofs
of the Chebotarev density theorem for number fields.
It has been used again in the proof of the Chebotarev density
theorem for function fields of one variable over finite fields,
of the arithmetic proof of the Hilbert irreducibility
theorem,
in the theory of Frobenius fields,
and on many more occasions.
See also [Deb99, Prop.~22.2].
Here we use it to replace the counting of points with a
given Artin class by the counting of rational points of an
absolutely irreducible variety over a finite field
(Theorem 3.1).
The argument itself appears in the proof of Lemma 2.8.

\reminfo {Definition 2.1}
{Ring-cover}
Let $S/R$ be an extension of integral domains whose
corresponding extension of quotient fields $F/E$ is finite and
separable.
We say that $S/R$ is a {\bf ring-cover}
if $R$ is integrally closed, $S=R[z]$, $z$ is integral over
$R$, and the discriminant of $\irr(z,E)$ is a unit of $R$.
In this case $S$ is the integral closure of $R$ in $F$
[FrJ08, Def.~1.6.3] and
$S/R$ is an \'etale extension of rings
[Ray70, p.~19].
Thus, the corresponding map
$\Spec(S)\to\Spec(R)$ of affine $K$-schemes is finite and
\'etale.
If $F/E$ is in addition Galois, we say that $S/R$ is a {\bf
Galois ring-cover}.
\endrem

For the rest of this section we fix a field $K$,
variables $A_1,\ldots,A_n,T$, let $\bfA=(A_1,\ldots,A_n)$,
and set $E=K(\bfA)$.

\rem {Definition 2.2}
We define the {\bf total degree} of a polynomial
$h=\sum_{i=1}^mc_iT^i\in E[T]$
to be
$$
\total.deg(h)=
\max_{1\le i\le m}\big(\max(\deg(f_i),\deg(g_i))+i\big),
$$
where $c_i={f_i(\bfA)\over g_i(\bfA)}$ is a reduced
presentation of $c_i$ with $f_i,g_i\in K[\bfA]$ and
$\deg(f_i),\deg(g_i)$ are the degrees of $f_i,g_i$,
respectively.
Under this definition, the total degree of $h$ coincides
with the degree of $h$ as a polynomial in $A_1,\ldots,A_n,T$
if $h\in K[\bfA,T]$.
\endrem

\reminfo {Definition 2.3}
{Bound}
In the following results we apply several algorithms to
polynomials $f_1,\ldots,f_k\in E[T]$ whose output are
polynomials $h_1,\ldots,h_{k'}\in E[T]$.
We say that the total degrees of $h_1,\ldots,h_{k'}$ are
{\bf bounded} if we can compute a function
$p\colon \bbZ\to\bbR$
that depends only on the algorithms (but not on $K$ neither
on $f_1,\ldots,f_k$) such that
$\sum_{i=1}^{k'}\total.deg(h_i)\le p(m)$,
where $m=\sum_{i=1}^k\total.deg(f_i)$.
Note that iteration of algorithms with output of bounded
total degree have again an output of bounded total degrees.
\endrem

\reminfo {Notation 2.4}
{Left conjugation}
Given a map $\tet$ from a set $G$ to a group $S$
and an element $\tau\in S$, we write
$^\tau\!\tet$ for the map from $G$ to $S$ defined for all
$\sig\in G$ by the rule
$^\tau\!\tet(\sig)=\tau\tet(\sig)\tau^{-1}$.
Note that $^{\tau\tau'}\!\tet={^\tau}\!(^{\tau'}\!\tet)$
for $\tau,\tau'\in S$ and $^1\tet=\tet$ for the unit element
$1$ of $S$.
\endrem

\reminfo {Definition 2.5}
{Points and homomorphisms}
Let $A$ be an integrally closed integral domain which is
finitely generated over a field $K$ and such that
$E=\Quot(A)$ is regular over $K$.
Let $F$ be a finite Galois extension of $E$ and write $B$
for the integral closure of $A$ in $F$.
Suppose $B/A$ is a ring cover.
Let $X=\Spec(A)$ and $Y=\Spec(B)$, and let $\lam\colon Y\to X$
be the finite \'etale morphism associated with the inclusion
$A\cont B$.
The points $\bfa\in X(K)$ bijectively corresponds
to $K$-homomorphisms
$\phi_\bfa\colon A\to K$.
The points $\bfb\in Y(K_s)$ with $\lam(\bfb)=\bfa$
bijectively correspond to epimorphisms $\phi_\bfb\colon B\to
K(\bfb)$ that extend $\phi_\bfa$.
Since $B/A$ is \'etale and Galois,
for each $\bfb$ as above the extension
$K(\bfb)/K$ is Galois that depends only on $\bfa$,
and $\phi_\bfb$ defines an embedding
$\phi^*\colon\Gal(K(\bfb)/K)\to\Gal(F/E)$ such that
$\phi_\bfb(\phi^*(\sig)x)=\sig(\phi_\bfb(x))$ for all
$\sig\in\Gal(K(\bfb)/K)$ and $x\in B$ [FrJ08, Lemma 6.1.4].
Let $\Gal(K)=\Gal(K_s/K)$ be the absolute Galois group of
$K$.
Composing $\phi^*$ with
$\res\colon\Gal(K)\to\Gal(K(\bfb)/K)$, we get a homomorphism
$\phi_\bfb^*\colon\Gal(K)\to\Gal(F/E)$
with $\Ker(\phi_\bfb^*)=\Gal(K(\bfb))$
such that
$\phi_\bfb(\phi_\bfb^*(\sig)x)=\sig(\phi_\bfb(x))$
for all $\sig\in\Gal(K)$ and $x\in B$.
The image $\phi_\bfb^*(\Gal(K))$ of $\phi_\bfb^*$ is
the {\bf decomposition group} $D_{\bfb/\bfa}$
of $\bfb$ over $\bfa$.
In particular, if $F'$ is the fixed field of $D_{\bfb/\bfa}$
in $F$ and $x\in B\cap F'$, then $\phi_\bfb(x)\in K$.
Finally we note that for each $\bfa\in X(K)$,
the action of $\Gal(F/E)$ on the set of prime ideals
of $B$ lying over $\Ker(\phi_\bfa)$ defines an action of
$\Gal(F/E)$ on $\lam^{-1}(\bfa)\cap Y(K_s)$
such that for all $\bfb\in Y(K_s)$ with $\lam(\bfb)=\bfa$
and each $\tau\in\Gal(F/E)$ we have
$K(\tau\bfb)=K(\bfb)$,
$\lam(\tau\bfb)=\bfa$,
$\phi_{\tau\bfb}=\phi_\bfb\circ\tau^{-1}$,
and $\phi_{\tau\bfb}^*={^\tau}\!(\phi_\bfb^*)$.

In order to prove the latter equality we consider
$\sig\in\Gal(K)$ and $x\in B$.
By definition
$$
\eqalign{
\phi_\bfb(\tau^{-1}(\phi_{\tau\bfb}^*(\sig)x))
&=\phi_{\tau b}(\phi_{\tau\bfb}^*(\sig)x)
=\sig(\phi_{\tau\bfb}(x))
\cr
&=\sig(\phi_\bfb(\tau^{-1}x))
=\phi_\bfb(\phi_\bfb^*(\sig)(\tau^{-1}x)).
}
$$
If $x$ is a primitive element of the cover $B/A$, then 
$\phi_\bfb$ maps the set of $K$-conjugates of $x$
injectively into $K_s$
(because the discriminant of $\irr(x,E)$ is a unit of $B$).
It follows in this case that
$\tau^{-1}(\phi_{\tau\bfb}^*(\sig)x)
=\phi_\bfb^*(\sig)(\tau^{-1}x)$.
Thus,
$\phi_{\tau\bfb}(\sig)x=(\tau\circ\phi_\bfb^*(\sig)\circ\tau^{-1})x$.
Therefore,
$\phi_{\tau\bfb}(\sig)x=\tau\circ\phi_\bfb^*(\sig)\circ\tau^{-1}$
for all $\sig\in\Gal(K)$.
This means that
$\phi_{\tau\bfb}(\sig)={^\tau}\!(\phi_\bfb^*)$,
as claimed.

The star operation is functorial in $B$:
Let $F'$ be a finite Galois extension of $E$ that contains
$F$,
let $B'$ be the integral closure of $A$ in $F'$,
and let $Y'=\Spec(B')$.
Suppose $B'/A$ is a ring-cover and let $\bfb'\in Y'(K_s)$
lie over $\bfb\in Y(K_s)$.
Then $\phi_{\bfb'}^*|_F=\phi_\bfb^*$.
\endrem

\proclaim {Lemma 2.6}
Let $t_1,\ldots,t_r$ be elements of $E_s$
and let $f_1,\ldots,f_r\in E[T]$ be separable polynomials
that satisfy
$f_i(\bfA,t_i)=0$, $i=1,\ldots,r$.
Then there exist polynomials $h\in K[\bfA,T]$ and
$0\ne g\in K[\bfA]$ of bounded total degrees
and an element $t\in E_s$
such that $h(\bfA,t)=0$, 
$E(t)=E(t_1,\ldots,t_r)$,
and $K[\bfA,g(\bfA)^{-1},t]/K[\bfA,g(\bfA)^{-1}]$ is a cover
of rings.
Moreover, the discriminant of $f_i(\bfA,T)$ as a polynomial
in $\bfA$ is invertible in $K[\bfA,g(\bfA)^{-1}]$,
$i=1,\ldots,r$.
\endproclaim

\demo {Proof}
The proof of the primitive element theorem [Lan93, Thm.~V.4.6]
assures that if
$C$ is a subset of $K[\bfA]$ of a large bounded cardinality,
then there exist $c_1,\ldots,c_r\in C$ such that
$t=c_1t_1+\cdots+c_rt_r$ satisfies
$K(t)=K(t_1,\ldots,t_r)$.
In particular, there exist such $c_1,\ldots,c_r\in K[\bfA]$
with bounded degrees.
Dividing each $f_i$ by its leading coefficient, we may
assume that $f_i$ is monic (as a polynomial in $T$).
Let
$$
h(\bfA,T)=\prod_{(t'_1,\ldots,t'_r)}(T-\sum_{i=1}^rc_it'_i),
$$
where for each $1\le i\le r$ the index $t'_i$ ranges over
all roots of $f_i$ in $E_s$.
Then $h(\bfA,T)\in E[T]$ and the coefficients of $h$ are
polynomials in the coefficients of $f_i$ with bounded
degrees.
Moreover, $h(\bfA,t)=0$.

Let $g(\bfA)$ be the product of the numerators and the
denominators of the discriminants of $f_1,\ldots,f_r,h$.
Then $g\in K[\bfA]$ and has a bounded degree.
By Definition 2.1, $h$ and $g$ satisfy all of the
requirements of the lemma.
\qed

\proclaim {Lemma 2.7}
Let $G$ be a group of order $d$,
let $\tet$ be the regular embedding of $G$ in $S_d$,
and let $H=\{\tau\in S_d\st ^\tau\!\tet=\tet\}$.
Then $|H|=d$.
\endproclaim

\demo {Proof}
We identify $S_d$ with the group $S_G$ of all permutations
of $G$.
Then $\tet(\sig)(x)=\sig x$ for all $\sig,x\in G$.
It follows that $\tau^{-1}\in H$ if and only if
$\sig\cdot\tau(x)=\tau(\sig x)$
(where the left hand side is the product of the elements
$\sig$ and $\tau(x)$ of $G$)
for all $\sig,x\in G$.
In particular, for $\sig=x^{-1}$ we get
$x^{-1}\cdot\tau(x)=\tau(1)$, so $\tau(x)=x\cdot\tau(1)$.
Conversely, if $\tau(x)=x\cdot\tau(1)$, then the former
condition holds, so $\tau\in H$.
Since $\tau(1)$ can take $d$ values, i.e.~the elements of
$G$, there are exactly $d$ possibilities for $\tau$.
\qed

\smallskip
The following central result is built on [BaS10, Lemma 2.2].

\proclaim {Lemma 2.8}
Let $K$ be a field and $A_1,\ldots,A_n,T$ variables.
For each $1\le i\le r$ let $f_i\in K[\bfA,T]$ be an
absolutely irreducible polynomial 
which is separable and of degree $d_i$ in $T$.
Let $L_i$ be a Galois
extension of $K$ of degree $d_i$.

We denote the splitting field of $f=f_1\cdots f_r$,
considered as a polynomial in $T$, over $E=K(\bfA)$ by $F$
and assume that $\Gal(F/E)\isom\prod_{i=1}^rS_{d_i}$.

Then there exist 
a proper algebraic subset $V$ of $\bbA_K^n$,
an absolutely irreducible normal affine
$K$-variety $W'$, and a finite \'etale map $\rho'\colon W'\to U$
with $U=\bbA_K^n\hefresh V$
such that the following conditions hold:

\condition (a)
$V$ and $W'$ are defined in $\bbA_K^n$ and $\bbA_K^{n+3}$,
respectively, by polynomials with coefficients in
$K$ of bounded total degrees.
\condition (b)
$\deg(f_i(\bfa,T))=d_i$ for each $\bfa\in U(K)$ and for $i=1,\ldots,r$.
\condition (c)
$\rho'(W'(K))$ is the set of all $\bfa\in U(K)$ such that
$L_i$ is generated by a root of $f_i(\bfa,T)$,
$i=1,\ldots,r$.
In particular, $f_i(\bfa,T)$ is irreducible of degree $d_i$ for
all $\bfa\in\rho'(W'(K))$ and $1\le i\le r$.
\condition (d)
$|(\rho')^{-1}(\bfa)\cap W'(K)|=\prod_{i=1}^rd_i$
for all $\bfa\in\rho'(W'(K))$.
\endproclaim

\demo {Proof}
For each $1\le i\le r$ let $F_i$ be the splitting
field of $f_i$ over $E$. Then $F=\prod_{i=1}^rF_i$ and since
$\Gal(F/E)\isom\prod_{i=1}^rS_{d_i}$, we have $\Gal(F_i/E)\isom
S_{d_i}$ for $i=1,\ldots,r$. For each $1\le i\le r$ we define a
homomorphism $\tet_i\colon\Gal(K)\to S_{d_i}$ in the following way:
First we identify the symmetric group $S_{d_i}$ with the group of
all permutations of $\Gal(L_i/K)$ (which we consider as a set of
$d_i$ elements). Then, for each $\sig\in\Gal(K)$ we define
$\tet_i(\sig)$ as a permutation of $\Gal(L_i/K)$ by the rule
$\tet_i(\sig)(\lam)=\sig|_{L_i}\cdot\lam$ for all
$\lam\in\Gal(L_i/K)$. Identifying $\Gal(F_i/E)$ with $S_{d_i}$, this
defines a homomorphism $\tet_i\colon\Gal(K)\to\Gal(F_i/E)$ with
$\Ker(\tet_i)=\Gal(L_i)$. Since
$\Gal(F/E)\isom\prod_{i=1}^r\Gal(F_i/E)$, the map
$\tet\colon\Gal(K)\to\Gal(F/E)$ defined by
$\tet(\sig)|_{F_i}=\tet_i(\sig)$ for all $\sig\in\Gal(K)$ and $1\le
i\le r$ is a well defined homomorphism with $\Ker(\tet)=\Gal(L)$,
where $L=L_1\cdots L_r$. It induces an embedding
$\tetgag\colon\Gal(L/K)\to\Gal(F/E)$ satisfying
$\tetgag(\sig|_L)=\tet(\sig)$ for each $\sig\in\Gal(K)$. We denote
the fixed field of $\tet(\Gal(K))$ in $F$ by $E'$.

Now we break up the rest of the proof into several parts.

\subdemoinfo {Part A}
{Field crossing argument}
We consider the Galois extension $\Fhat=FL$ of $E$ and
note that since $F/K$ is a regular extension,
$$
\Gal(\Fhat/E')=\Gal(\Fhat/F)\times\Gal(\Fhat/E'L)
\isom\Gal(L/K)\times\Gal(F/E'),
\leqno (\cf1)
$$
where the latter isomorphism is induced by the corresponding
restriction maps.
Next we consider the subgroup
$$
\Del=\{\del\in\Gal(\Fhat/E)\st\tetgag(\del|_L)=\del|_F\}
$$
of $\Gal(\Fhat/E')$.
Using (\cf1) and the injectivity of $\tetgag$,
one sees that
$\Del\cap\Gal(\Fhat/F)=1$
and
$\Del\cap\Gal(\Fhat/E'L)=1$.
Moreover,
$\Del\cdot\Gal(\Fhat/F)=\Gal(\Fhat/E)$ and
$\Del\cdot\Gal(\Fhat/E'L)=\Gal(\Fhat/E')$.

Indeed, by (\cf1),
for each $\eps\in\Gal(\Fhat/E')$
there exists $\del\in\Gal(\Fhat/E')$ such that
$\del|_F=\eps|_F$ and $\del|_L=(\tetgag)^{-1}(\eps|_F)$.
Then, $\eps=\del\cdot\del^{-1}\eps$,
\ $\del\in\Del$, and $(\del^{-1}\eps)|_F=1$, so
$\del^{-1}\eps\in\Gal(\Fhat/F)$.
Similarly,
there exists
$\del'\in\Gal(\Fhat/E')$ such that $\del'|_L=\eps|_L$
and $\del'|_F=\tetgag(\eps|_L)$.
Thus, $\eps=\del'\cdot(\del')^{-1}\eps$, \
$\del'\in\Del$, and
$((\del')^{-1}\eps)|_L=1$, so
$(\del')^{-1}\eps\in\Gal(\Fhat/E'L)$.
It follows that the fixed field $F'$ of $\Del$ in $\Fhat$
satisfies
$$
FF'=F'L=\Fhat
\qquad\hbox{and}\qquad
F\cap F'=F'\cap E'L=E'
\leqno (\cf2)
$$
and fits into the following diagram of fields:
$$
\xymatrix@R=10pt{
F \ar@{-}[rr]
&& \Fhat
\\
& F' \ar@{-}[ur]^\Del
\\
E' \ar@{-}[uu] \ar@{-}[ur] \ar@{-}[rr]
&& E'L \ar@{-}[uu]
\\
E \ar@{-}[u] \ar@{-}[rr]
&& EL \ar@{-}[u]
\\
K \ar@{-}[u] \ar@{-}[rr]
&& L \ar@{-}[u]
}
$$

\subdemoinfo {Part B}
{Integral \'etale extensions of rings}
By Lemma 2.6 applied to the roots of $f$,
there exist a
polynomial $h\in K[\bfA,T]$, separable in $T$, and a nonzero
polynomial $g\in K[\bfA]$ of bounded total degrees,
and there exists an element $t\in F$ such that
$h(\bfA,t)=0$, $F=E(t)$,
and
$K[\bfA,g(\bfA)^{-1},t]/K[\bfA,g(\bfA)^{-1}]$
is a ring-cover.
Moreover, the discriminant of each $f_i(\bfA,T)$ considered
as a polynomial in $T$ is invertible in
$K[\bfA,g(\bfA)^{-1}]$.

Applying Lemma 2.6
to $t$ and to a primitive element of $L/K$,
we find a polynomial $\hhat\in K[\bfA,T]$, separable in $T$,
a nonzero polynomial $\ghat\in K[\bfA]$, 
and an element
$\that\in\Fhat$ such that
$\hhat,\ghat$ have bounded degrees,
$\hhat(\bfA,\that)=0$,
$\Fhat=E(\that)$, and
$K[\bfA,\ghat(\bfA)^{-1},\that]
/K[\bfA,\ghat(\bfA)^{-1}]$ is a
ring-cover.
Moreover, $\ghat$ may be taken as a multiple of $g$.

Let $\that_1,\ldots,\that_k$ be a $\Del$-orbit starting with
$\that_1=\that$ and write
$\prod_{j=1}^k(T-\that_j)=T^k+t'_1T^{k-1}+\cdots+t'_k$.
Then $F'=E(t'_1,\ldots,t'_k)$ [FrJ08, End of the proof of
Lemma 19.3.2].
Moreover, each of the coefficients of $\irr(t'_j,E)$,
$j=1,\ldots,k$,
is a symmetric polynomial in the
$E$-conjugates of $\that$ of a bounded total degree.
Hence, those coefficients are polynomials in the coefficients
of $\hhat$ of bounded degrees having integral coefficients.
Applying Lemma 2.6
now to $\irr(t'_1,E),\ldots,\irr(t'_k,E)$,
we can compute polynomials
$h'\in K[\bfA,T]$ and $g'\in K[\bfA]$ of bounded total degrees
and an element $t'\in F'$
such that $h'(\bfA,t')=0$, $g'(\bfA)\ne0$, \
$S'=K[\bfA,g'(\bfA)^{-1},t']$
is a ring cover of
$R=K[\bfA,g'(\bfA)^{-1}]$,
and $F'=E(t')$.
Moreover, $g'$ may be taken as a multiple of $\ghat$, hence
of $g$.
Also,
\condition (\cf3)
the discriminant of $f_i(\bfA,T)$ is invertible in $R$,
$i=1,\ldots,r$.

Let $S=R[t]$ and $\Shat=R[\that]$.
Then $S$ is the integral closure of $R$ in $F$,
$S'$ is the integral closure of $R$ in $F'$,
and $\Shat$ is the integral closure of $R$, $S$, and $S'$ in
$\Fhat$.
It follows from (\cf2) that $SL=\Shat$ and $S'L=\Shat$.
Moreover, $S/R$, $S'/R$,
$\Shat/S$, $\Shat/S'$, and $\Shat/R$ are ring covers in the sense of
Definition 2.1.

\subdemoinfo {Part C}
{Fiber products}
Let $V$ be the Zariski closed subset of $\bbA_K^n$
defined by the equation $g'(\bfA)=0$.
Then $U=\bbA_K^n\hefresh V$ is a non-empty Zariski open
subset of $\bbA_K^n$.
We set $W=\Spec(S)$, $W'=\Spec(S')$, and
$\What=\Spec(\Shat)$.
By Part B, $V$, $W$, and $W'$
are closed subsets of $\bbA_K^n$, $\bbA_K^{n+3}$,
and $\bbA_K^{n+3}$, respectively, defined
by polynomials with coefficients in $K$
of bounded degrees.
Moreover, the ring-covers in the left diagram of
inclusions of rings in (\cf4) define finite \'etale
morphisms as in the right diagram of (\cf4):
$$
\xymatrix{
S \ar@{-}[rr]
&& \Shat
\\
& S' \ar@{-}[ur]
\\
R \ar@{-}[uu] \ar@{-}[ur] \ar@{-}[rr]
&& RL \ar@{-}[uu]
}
\hskip 2cm
\xymatrix{
W \ar[dd]_\rho
&& \What \ar[ll]_\ome \ar[dl]_{\ome'} \ar[dd]^{\rho_L}
\\
& W' \ar[dl]_{\rho'}
\\
U
&&
U_L \ar[ll]
}
\leqno (\cf4)
$$
Here $U_L=U\times_{\Spec(K)}\Spec(L)=\Spec(RL)$ and
each of the three rectangles in the right diagram of (\cf4)
is cartesian, i.e.~$\What=W\times_UW'$,
$\What=W\times_UU_L$,
and $\What=W'\times_UU_L$.
This concludes the proof of (a).

\subdemoinfo {Part D}
{Geometric points}
For each $\bfa\in U(K)$
and each $\bfb\in W(K_s)$ with $\rho(\bfb)=\bfa$
we consider the homomorphisms
$\phi_\bfb\colon S\to K_s$ and
$\phi_\bfb^*\colon\Gal(K)\to\Gal(F/E)$ introduced in Definition
2.5.
By Diagram (\cf4), $\phi_\bfb\colon S\to K_s$
uniquely extends to an $L$-homomorphism
$\phihat\colon\Shat\to K_s$.
Let $\bhat$ be the unique point of $\What(K_s)$ with
$\ome(\bfbhat)=\bfb$ and
$\phi_\bhat=\phihat$.
Then let $\bfb'$ be the unique point of $W'(K_s)$ with
$\ome'(\bfbhat)=\bfb'$ and
$\phi_{\bfb'}=\phi_\bfbhat|_{S'}$.
Again, the map $\bfbhat\onto\bfb'$ is bijective.

\subdemo {Claim}
$\bfb'\in W'(K)$ if and only if $\phi_\bfb^*=\tet$.

First we note that by (\cf2),
there exists for each $\sig\in\Gal(K)$ a unique element
$\tethat(\sig)\in\Del=\Gal(\Fhat/F')$ such that
$$
\tethat(\sig)|_F=\tet(\sig)
\hbox { and }
\tethat(\sig)|_L=\sig|_L.
\leqno (\cf5)
$$
Since $\tet(\Gal(K))=\Gal(F/F')$, this implies that
$\tethat\colon\Gal(K)\to\Gal(\Fhat/F')$
is an epimorphism with $\Ker(\tethat)=\Gal(L)$.

Suppose $\bfb'\in W'(K)$.
Then $\bfbhat\in\What(L)$, so
$\phi_\bfbhat$ fixes the elements of $L$.
Hence, for each $\sig\in\Gal(K)$ and $x\in L$ we have
$\phi_\bfbhat^*(\sig)x
=\phi_\bfbhat(\phi_\bfbhat^*(\sig)x)
=\phi_\bfbhat(\sig x)
=\sig x$.
Therefore,
$\phi_\bfbhat^*(\sig)|_L=\sig|_L=\tethat(\sig)|_L$.
Also,
$\phi_\bhat^*(\Gal(K))\le\Gal(F/F')$, so that
$\phi_\bfbhat^*(\sig)x
=x
=\tethat(\sig)x$ for each $x\in F'$.
It follows from (\cf1) that $\phi_\bfbhat^*=\tethat$.
Applying $\res_{\Fhat/F}$ to the latter equality, we get
$\phi_\bfb^*=\tet$.

Conversely, suppose $\phi_\bfb^*=\tet$.
Then $\Ker(\phi_\bfb^*)=\Ker(\tet)=\Gal(L)$,
so $\bfb\in W(L)$, hence $\bfbhat\in\What(L)$.
It follows that for each $\sig\in\Gal(K)$ we have
$\phi_\bhat^*(\sig)|_F=\phi_\bfb^*(\sig)=\tet(\sig)
=\tethat(\sig)|_F$
and
$\phi_\bhat^*(\sig)|_L
=\sig|_L
=\tethat(\sig)|_L$.
By Part A,
$\phi_\bhat^*(\sig)=\tethat(\sig)\in\Gal(\Fhat/F')$.
Therefore, $\phi_\bhat(S')=K$.
Consequently, $\bfb'=\ome'(\bfbhat)$ satisfies $\bfb'\in W'(K)$.
This concludes the proof of the Claim.

\subdemoinfo {Part E}
{The orbit of $\tet$}
Following Notation 2.4,
the group $\Gal(F/E)$ acts on the set of all homomorphisms
from $\Gal(K)$ into $\Gal(F/E)$ by left conjugation.
Let $\Tet=\{^\tau\!\tet\st\tau\in\Gal(F/E)\}$ be the
$\Gal(F/E)$-orbit of $\tet$.

For each $\bfa\in U(K)$ 
we consider 
the conjugacy class of homomorphisms
$\Phi_\bfa^*=\{\phi_\bfb^*\st\bfb\in W(K_s)\hbox{ and }
\rho(\bfb)=\bfa\}$. By Part D
$$
\rho'(W'(K))=\{\bfa\in U(K)\st\Phi_\bfa^*=\Tet\}.
\leqno (\cf6)
$$

\subdemoinfo {Claim}
{A point $\bfa\in U(K)$ satisfies $\Phi_\bfa^*=\Tet$ if and
only if each root of $f_i(\bfa,T)$ generates $L_i$ over $K$,
$i=1,\ldots,r$}

First suppose $\Phi_\bfa^*=\Tet$.
Then there exists $\bfb\in W(K_s)$ such that
$\rho(\bfb)=\bfa$ and
$\phi_\bfb^*=\tet$.
Since the discriminant of $f_i(\bfA,T)$ is a unit of $R$,
\ $\phi_\bfb$ maps the set of roots of $f_i(\bfA,T)$ in $F$
bijectively on the set of roots of $f_i(\bfa,T)$ in $K_s$.
Let $x$ be a root of $f_i(\bfA,T)$ in $F$ and let
$\sig\in\Gal(K)$.
If $\sig(\phi_\bfb(x))=\phi_\bfb(x)$,
then $\phi_\bfb(\phi_\bfb^*(\sig)x)=\phi_\bfb(x)$,
hence, by the injectivity of $\phi_\bfb$, we have
$\phi_\bfb^*(\sig)x=x$, so $\tet(\sig)x=x$.
By the definition of $\tet$, there exists $\lam\in\Gal(L_i/K)$
such that $\sig|_{L_i}\lam=\lam$.
Therefore, $\sig\in\Gal(L_i)$.
Conversely, if $\sig\in\Gal(L_i)$, then
$\sig(\phi_\bfb(x))
=\phi_\bfb(\phi_\bfb^*(\sig)x)
=\phi_\bfb(\tet(\sig)x)
=\phi_\bfb(\tet_i(\sig)x)
=\phi_\bfb(x)$.
Consequently, $L_i=K(\phi_\bfb(x))$.

Conversely, suppose for each $1\le i\le r$ there exists
a root $\xgag_i$ of $f_i(\bfa,T)$ with
$L_i=K(\xgag_i)$.
Then, by (b) and by the assumption on $L_i$,
$\deg(f_i(\bfa,T)=d_i=[L_i:K]=[K(\xgag_i):K]$.
Hence, $\Gal(L_i/K)$ acts freely on the roots of $f_i(\bfa,T)$.
Choosing
$\bfb\in W(K_s)$ with $\rho(\bfb)=\bfa$,
this implies that the group $\phi_\bfb^*(\Gal(K))|_{F_i}$
acts freely on the roots of
$f_i(\bfA,T)$.
By definition, this is also the case for the group
$\tet(\Gal(K))|_{F_i}$.
Since
 $\Gal(f_i(\bfA,T),K(\bfA))=S_{d_i}$, there exists
$\tau_i\in\Gal(F_i)$ such that
$^{\tau_i}\!(\res_{F/F_i}\circ\tet)
=\res_{F_i/F}\circ\phi_\bfb^*$.
Since $\Gal(F/E)=\prod_{i=1}^r\Gal(F_i/E)$, there exists
$\tau\in\Gal(F/E)$ such that $\tau|_{F_i}=\tau_i$ for
$i=1,\ldots,r$.
Hence, $^\tau\!\tet=\phi_\bfb^*$.
Therefore, $\Tet=\Phi_\bfa^*$, as claimed.

The first statement of (c) of our Lemma
follows now from (\cf6) and from the Claim.
The second one follows from the first and from (b).

\subdemoinfo {Part F}
{Let $\bfa\in U(K)$.
We prove that
$|(\rho')^{-1}(\bfa)\cap W'(K)|={[F:E]\over|\Tet|}$}

The stabilizer $H=\{\tau\in\Gal(F/E)\st^\tau\!\tet=\tet\}$
of $\tet$ satisfies
$$
|H|={[F:E]\over|\Tet|}
\leqno (\cf7)
$$

Let $B=\{\bfb\in\rho^{-1}(\bfa)\st\phi_\bfb^*=\tet\}$.
We prove that $H$ acts regularly on $B$.
Indeed, if $\tau\in H$ and $\bfb\in B$, then by Definition
2.5,
$\phi_{\tau\bfb}^*={^\tau\!(\phi_\bfb^*)}={^\tau\!\tet}=\tet$,
so $\tau\bfb\in B$.
Thus, $H$ acts on $B$.
Next we prove that the action is transitive.
To this end let $\bfb,\bfb'\in B$.
Then there exists $\tau\in\Gal(F/E)$ with
$\bfb'=\tau\bfb$.
Hence, $\tet=\phi_{\bfb'}^*=\phi_{\tau\bfb}^*
={^\tau\!(\phi_\bfb^*)}
={^\tau\!\tet}$, so $\tau\in H$, as desired.
Finally the action of $H$ is free.
Indeed, if $\bfb\in B$ and $\tau\bfb=\bfb$, then by
Definition 2.5,
$\phi_\bfb\circ\tau^{-1}=\phi_{\tau\bfb}=\phi_\bfb$.
Since $\phi_\bfb$ is injective on the roots of $h(\bfA,T)$,
this implies that
$\tau=1$, which proves our assertion.

It follows that $|B|=|H|$.
As in Part D let $B'$ be the set of all points $\bfb'\in
W'(K_s)$ corresponding to the points $\bfb\in B$.
Since the map $\bfb\onto\bfb'$ is injective, we have
$|B'|=|B|=|H|$.
Moreover, by the Claim of Part D,
$B'=\rho^{-1}(\bfa)\cap W'(K)$.
Hence, by (\cf7),
$|\rho^{-1}(\bfa)\cap W'(K)|={[F:E]\over|\Tet|}$,
as claimed.

\subdemoinfo {Part G}
{We prove:
$|\Tet|=\prod_{i=1}^r(d_i-1)!$}

By (\cf7),
$|\Tet|={|\Gal(F/E)|\over|H|}$.
Since
$|\Gal(F/E)|=|\prod_{i=1}^rS_{d_i}|=\prod_{i=1}^rd_i!$,
it suffices to prove that $|H|=\prod_{i=1}^rd_i$.
Since
$\Gal(F/E)\isom\prod_{i=1}^r\Gal(F_i/E)$
and $\tet=\prod_{i=1}^r\tet_i$, we have
$^\tau\!\tet=\tet$ if and only if $^{\tau_i}\tet_i=\tet_i$,
where $\tau_i=\tau|_{F_i}$ for $i=1,\ldots,r$.
Hence, it suffices to prove that
$\#\{\tau\in\Gal(F_i/E)\st^\tau\!\tet_i=\tet_i\}=d_i$.
But this follows from Lemma 2.7.

The combination of Part F and Part G implies Statement (d)
of our lemma.
\qed

\section{Finite Fields and PAC Fields}

We combine Corollary 1.3 with Lemma
2.8 to the case when $K=\bbF_q$ 
and establish an asymptotic formula for the
number of pairs $(a_1,a_2)\in\bbF_q^2$ such that
$f_i(T,a_1T+a_2)$ is irreducible in $\bbF_q[T]$ for
$i=1,\ldots,r$,
when $f_1,\ldots,f_r\in\bbF_q[T,X]$ are
characteristic-$0$-like nodal non-associate
polynomials.

\proclaim {Theorem 3.1}
Let 
$f_1,\ldots,f_r\in\bbF_q[T,X]$ be characteristic-$0$-like
nodal non-associate polynomials, where
$q$ is a power of a prime $p$ such that $r=1$ if $p=2$. 
For each $1\le i\le r$ let $d_i=\deg(f_i(T,X))$
and set $d=\prod_{i=1}^rd_i$.
Then,
$$
\eqalign{
\#\{(a_1,a_2)\in\bbF_q^2\st f_i(T,a_1T+a_2)
\hbox{ is irreducible in } \bbF_q[T],
\ &i=1,\ldots,r\}
\cr
&={q^2\over d}+O(q^{3/2}),
}
$$
where the constant of the $O$ is a computable function in
$\sum_{i=1}^r\deg(f_i(T,X))$.
\endproclaim

\demo {Proof}
Let $K=\bbF_q$ and
let $A_1,A_2$ be additional variables.
For each $1\le i\le r$ let
$f'_i(A_1,A_2,T)=f_i(T,A_1T+A_2)$,
let $F_i$ be the the splitting field
of $f'_i(A_1,A_2,T)$ over $E=K(A_1,A_2)$, and let $F=F_1\cdots F_r$.
By Corollary 1.3,
$\Gal(F/E)=\prod_{i=1}^rS_{d_i}$.
For each $1\le i\le r$ 
let $L_i=\bbF_{q^{d_i}}$ be the unique Galois extension
of $\bbF_q$ of degree $d_i$ and note that
$d_i=\deg_T(f'_i(A_1,A_2,T))$.

Let
$m=\sum_{i=1}^r\deg(f'_i(A_1,A_2,T))=\sum_{i=1}^r\deg(f_i(T,X))$
and consider the objects $V$, $U=\bbA_K^2\hefresh V$,
and $\rho'\colon W'\to U$ 
supplied by Lemma 2.8 with respect to
$K=\bbF_q$, to $f'_1,\ldots,f'_r$, and to $L_1,\ldots,L_r$.
In particular, 
$V$ is a Zariski closed subset of $\bbA_K^2$ defined by a
nonzero polynomial of an $m$-bounded degree.
Thus, there exists an $m$-bounded constant $c_1$ such that
$$
\#V(K)\le c_1q.
\leqno (\cf1)
$$
Moreover, $W'$ is an absolutely irreducible affine 
$K$-subvariety of $\bbA_K^5$
 defined by polynomials of
$m$-bounded degrees.
In addition,
$\dim(W')=2$, because $\rho'\colon W'\to U$ is a
finite morphism.
The Lang-Weil estimates give an explicit $m$-bounded
constant $c_2$
such that
$$
|\#W'(K)-q^2|\le c_2q^{3\over2}.
\leqno (\cf2)
$$
(See [LaW54, Thm.~1], [FHJ94, (1) and (2) of Section 3],
or [Zyw10, Thm.~2.1]).
Note that the estimates given in [Zyw10] are exponential in
the degrees of the polynomials that define $W'$.)

Let $I$ be the set of all $(a_1,a_2)\in\bbF_q^2$ such that
$f_i(T,a_1T+a_2)$ is irreducible in $\bbF_q[T]$ for
$i=1,\ldots,r$.
By Lemma 2.8(c),
$\rho'(W'(K))$ is the set of all $(a_1,a_2)\in U(K)$ such
that for each $1\le i\le r$
the field $L_i$ is generated over $K$ by a root of $f_i(T,a_1T+a_2)$.
Moreover, by Lemma 2.8(b),
$\deg(f_i(T,a_1T+a_2))=d_i$ for $i=1,\ldots,r$.
Since $L_i$ is the unique extension of $K$ of degree $d_i$,
this implies that $\rho'(W'(K))=I\cap U(K)$,
so $|\rho(W'(K))|=d\cdot|(I\cap U(K))|$.

By Lemma 2.8(d),
$|(\rho')^{-1}(\bfa)\cap W'(K)|=d$ for each
$\bfa\in\rho(W'(K))$.
Hence, by (\cf2),
$|\#(I\cap U(K))-{q^2\over d}|\le{c_2\over d}q^{3\over2}$.
It follows from (\cf1) that
$|\#I-{q^2\over d}|\le {c_2\over d}q^{3\over2}+c_1q$,
as desired.
\qed

\smallskip
Almost the same proof can be applied to PAC fields.

\proclaim {Theorem 3.2}
Let $K$ be a PAC field
and let $f_1,\ldots,f_r\in K[T,X]$ be
charac\-teristic-$0$-like nodal polynomials.
Suppose for each $i$ the field $K$ has a Galois extension
$L_i$ of degree $\deg(f_i(T,X))$.
Then $\bbA_K^2$ has a Zariski-dense subset $B$ such that
$f_i(T,a_1T+a_2)$ is irreducible for $i=1,\ldots,r$.
\endproclaim

\demo {Proof}
If one of two associated polynomial is irreducible, so is
the other.
Thus, after possible dropping some of the $f_i$'s, we may
assume that $f_1,\ldots,f_r$ are non-associate.
Let $A_1,A_2$ be additional variables and set $E=K(A_1,A_2)$.
For each $1\le i\le r$ let
$f'_i(A_1,A_2,T)=f_i(T,A_1T+A_2)$,
let $F_i$ be the splitting field of $f'_i(A_1,A_2,T)$ over
$E=K(A_1,A_2)$ and set $F=F_1\cdots F_r$.
By Corollary 1.3, 
$\Gal(F/E)=\prod_{i=1}^rS_{d_i}$ and note that
$d_i=\deg_T(f'_i(A_1,A_2,T))$.

Consider the objects $U$ and $\rho'\colon W'\to U$ supplied
by Lemma 2.8 with respect to
$f'_i,\ldots,f'_r$.
Since $K$ is PAC and $W'$ is an absolutely irreducible
$K$-variety, the set $W'(K)$ is non-empty,
so there exists $\bfa\in\rho'(W'(K))$.
By Lemma 2.8, $f_i(T,a_1T+a_2)$ is irreducible
for each $1\le i\le r$.
\qed

\section{Concluding Remarks}

The following proposition proves that Theorem 1.1 of
[BeW05] is a special case of our main result.

\proclaim {Proposition 4.1}
Let $K$ be a field of characteristic $p$
and $f\in K[T,X]$ an irreducible
polynomial of degree $d$ such that
\condition (\cf1a)
$p\notdivide d(d-1)$, and
\condition (\cf1b)
the Zariski closure $\Gam_\Kgal$
in $\bbP_\Kgal^2$ of the affine plane curve defined
by $f(T,X)=0$ is smooth.

\noindent
Then $\Gam$ is characteristic-$0$-like and nodal.
\endproclaim

\demo {Proof}
By assumption, $\Gam$ is irreducible and by (\cf1b), $\Gam$
is absolutely irreducible.
In the proof of [BeW05, Prop.~3.1], Bender and Wittenberg
show for $K=\bbF_q$
that assumption (\cf1) implies that the map of $\Gam_\Kgal$ into
its dual curve is separable.
The proof is however valid for every field $K$.
It follows from [Kat73, Cor.~3.5.0 and Cor.~3.2.1] that
the intersection multiplicities
of all but finitely many lines $L$ in $\bbP_\Kgal^2$ 
with $\Gam$ are at most $2$.
In particular, for only finitely many points
$\bfp\in\Gam_\Kgal$, the intersection number of the tangent
to $\Gam$ at $\bfp$ is greater than $2$.
This means that $\Gam$ has only finitely many inflection
points.

By (\cf1a), $p\ne2$.
Hence, by [GeJ89, Prop.~4.5], $\Gam_\Kgal$ has only finitely
many double tangents.
Again, by (\cf1a), $\Gam_\Kgag$ is not a line and not a
conic in characteristic $2$.
Hence, by Samuel [Har77, Thm.~IV.3.9], $\Gam_\Kgal$ is not
strange.
Finally, $\Gam_\Kgal$ is nodal because it is smooth.
Consequently, $\Gam$ is characteristic-$0$-like and nodal.
\qed

\smallskip
One of the ingredients of the proof of Lemma
1.1 (on which eventually our proof of Theorem
3.1 in case $r\ge2$ relies) is that the dual curves
of distinct absolutely irreducible projective plane
$K$-curves are distinct, if
$\chr(K)\ne2$.
The following example shows that this is not the case if
$\chr(K)=2$.
Thus, our proof of Theorem 3.1 fails if
$\chr(K)=2$ and $r\ge2$.
We do not know if the theorem itself holds in that case.

\rem {Example 4.2 (Bjorn Poonen)}
Let $K$ be an algebraically closed field of characteristic
$2$.
We consider the homogeneous polynomial
$$
f(X_0,X_1,X_2)=X_1^5+X_1^2X_0^3+X_0^4X_2
$$
and the projective plane curve $\Gam$ defined by the
equation $f(X_0,X_1,X_2)=0$.
The finite part of $\Gam$ is defined by the equation
$X^5+X^2=Y$.
This may be used to prove that $\Gam$ is absolutely
irreducible.
Moreover, $\Gam$ has only finitely many inflection
points (e.g.~$(0,0)$ is a simple non-inflection point, now use
[GeJ89, Cor.~3.2]),
only finitely many double tangents (see next paragraph),
and $\Gam$ is not
strange (the latter assertion is also a special case of a
theorem of Samuel [Har77, Thm.~IV.3.9]).

Let $(x_0{:}x_1{:}x_2)$ be a generic homogeneous point of
$\Gam$.
Applying (1) of Section 1, we find that 
$(x_0^*{:}x_1^*{:}x_2^*)=(x_0^2x_1^2:x_1^4:x_0^4)$
is a generic point of $\Gam^*$.
In particular the map $\Gam\to\Gam^*$ is purely
inseparable, so $\Gam$ has only finitely many double
tangents [GeJ89, Lemma 4.2].
The point $(x_0^*{:}x_1^*{:}x_2^*)$
lies on the irreducible homogeneous plane curve defined
by the equation $X_1X_2=X_0^2$.
Hence, the latter equation defines $\Gam^*$.

Next we define an additional projective pland curve $\Del$
by exchanging the variables $X_1$ and $X_2$.
In other words, the equation that defines $\Del$ is
$X_2^5+X_2^2X_0^3+X_0^4X_1=0$.
The equation that defines $\Del^*$ will therefore be
$X_2X_1=X_0^2$.
It follows that $\Gam^*=\Del^*$ although $\Gam\ne\Del$.
\endrem


\references {AAAA}

\ref [BaH62]
P. T. Bateman and R. A. Horn,
\sl A heuristic asymptotic formula concerning the
distribution of prime numbers,
\rm Mathematics of Computation {\bf 16} (1962), 363--367.

\ref [BaS10]
L. Bary-Soroker,
\sl Irreducible values of polynomials,
\rm manuscript 2010.

\ref [BaD04]
P. T. Bateman and H. G. Diamond,
\sl Analytic Number Theory - An Introductory Course,
\rm World Scientific Publishing Co., Denver 2004

\ref [BeW05]
A. O. Bender and O. Wittenberg
\sl A Potential analogue of Schinzel's hypothesis for
polynomials with coefficients in $\bbF_q[t]$,
\rm International Mathematics Research Notices {\bf 36}
(2005), 2238--2248.

\ref [CCG08]
B. Conrad, K. Conrad, and R. Gross,
\rm Prime specialization in genus $0$,
\rm Transactions of the AMS {\bf 360} (2008), 2867-2908.

\ref [Deb99]
P. D\`ebes,
\sl Galois Covers with Prescribed Fibers: The Beckmann-Black
Problem,
\rm Annale Scuola Normale Superiore Pisa {\bf 28} (1999),
273--286.

\ref [FHJ94]
M. Fried, D. Haran, M. Jarden,
\sl Effective counting of the points of definable sets over
finite fields,
\rm Israel Journal of Mathematics {\bf 85} (1994), 103--133.

\ref [FrJ76] 
M. Fried and M. Jarden,
\sl Diophantine properties of subfields of  $\bbQgal$, 
\rm American Journal of Mathematics {\bf 100} (1978), 653--666.

\ref [FrJ08]
M. D. Fried and M. Jarden,
\sl Field Arithmetic,
Third Edition, revised by Moshe Jarden,
\rm Ergebnisse der Mathematik (3) {\bf 11}, Springer,
Heidelberg, 2008.

\ref [Gau03]
S. Gao,
\sl Problem 6,
AIM Workshop on ``Future Directions in Algorithmic Number
Theory,''
\rm American Institute of Mathematics, California,
2003,\break
http://www.aimath.org/ARCC/workshops/primesinp.html

\ref [GeJ89]
W.-D.\ Geyer and M. Jarden,
\sl On stable fields in positive characteristic,
\rm Geometriae Dedicata {\bf 29} (1989), 335--375.

\ref [Ful89]
W. Fulton,
\sl Algebraic Curves,
\rm Addison Weseley, Redwood City, 1989.

\ref [Har77]
R. Hartshorne,
\sl Algebraic Geometry,
\rm Graduate Texts in Mathematics {\bf 52}, Springer, New York,
1977.

\ref [Kat73]
N. Katz,
\sl Pinceaux de Lefschetz: th\'eor\`eme d'existence, Groups
de monodromies en G\'eom\'etrie Alg\'ebrique II,
\rm Springer Lecture Notes in Mathematics {\bf 340} (1973),
212-353.

\ref [Lan93]
S. Lang,
\sl Algebra, Third Edition,
\rm Eddison-Wesley, Reading, 1993.

\ref [LaW54]
S. Lang and A. Weil,
\sl Number of points of varieties in finite fields,
\rm American Journal of Mathematics {\bf 76} (1954), 819--827.

\ref [Pol08]
P. Pollack,
\sl Simultaneous prime specializations of polynomials over
finite fields,
\rm Proceedings of the London Mathematical Society {\bf 97}
(2008), 545--567.

\ref [Ray70]
M. Raynaud,
\sl Anneaux Locaux Hens\'eliens,
\rm Lecture Notes in Mathematics {\bf 169}, Springer,
Berlin, 1970.

\ref [Swa62]
R. G. Swan,
\sl Factorization of polynomials over finite fields,
\rm Pacific Journal of Mathematics {\bf 12} (1962),
1099--1106.


\end{document}